\documentclass[11pt]{amsart}
\usepackage{amsmath,amssymb,latexsym,cite}
\usepackage[small]{caption}
\usepackage{graphicx,color,mathrsfs}%wasysym,overpic,tikz,
\usepackage{subfigure,color}
\usepackage{cite}
\usepackage[colorlinks=true,urlcolor=blue,
citecolor=red,linkcolor=blue,linktocpage,pdfpagelabels,
bookmarksnumbered,bookmarksopen]{hyperref}
\usepackage[italian,english]{babel}
\usepackage{units}
\usepackage{enumitem}
\usepackage[left=2.75cm,right=2.75cm,top=3.2cm,bottom=3.2cm]{geometry}

\numberwithin{equation}{section}
\newtheorem{theorem}{Theorem}[section]
\newtheorem{proposition}[theorem]{Proposition}
\newtheorem{lemma}[theorem]{Lemma}
\newtheorem{remark}[theorem]{Remark}

\newtheorem{definition}[theorem]{Definition}
\theoremstyle{definition}

\renewcommand{\epsilon}{\eps}
\renewcommand{\i}{{\rm i}}

\newcommand{\R}{{\mathbb R}}

\newcommand{\eps}{\varepsilon}

\newcommand{\Ma}{\mathscr{H}_{A,V}}

\newcommand{\pnorm}[2][]{\if #1'' \left|#2\right|_p \else \left|#2\right|_{#1} \fi}

\renewcommand{\theta}{\vartheta}

\title[Nonlocal Schr\"{o}dinger--Kirchhoff equations with magnetic field]{Nonlocal Schr\"{o}dinger--Kirchhoff equations\\with external magnetic field}
\author[X.\ Mingqi]{Xiang Mingqi}
\address[X.\ Mingqi]{College of Science
\newline\indent
Civil Aviation University of China
\newline\indent
Tianjin, 300300, P.R.\ China}
\email{\href{mailto:xiangmingqi\_hit@163.com}{xiangmingqi\_hit@163.com}}

\author[P.\ Pucci]{Patrizia Pucci}

\address[P.\ Pucci]{Dipartimento di Matematica e Informatica
\newline\indent
Universit\`a degli Studi di Perugia
\newline\indent
Via Vanvitelli 1, I-06123 Perugia, Italy}
\email{\href{mailto:patrizia.pucci@unipg.it}{patrizia.pucci@unipg.it}}

\author[M.\ Squassina]{Marco Squassina}

\address[M.\ Squassina]{Dipartimento di Informatica \newline\indent
Universit\`a degli Studi di Verona
\newline\indent
Strada Le Grazie 15, I-37134 Verona, Italy}
\email{\href{mailto:marco.squassina@univr.it}{marco.squassina@univr.it}}

\author[B.\ Zhang]{Binlin Zhang}

\address[B.\ Zhang]{Department of Mathematics \newline\indent
	Heilongjiang Institute of Technology
	\newline\indent
	Harbin, 150050, P.R.\ China}
\email{\href{mailto:zhangbinlin2012@163.com}{zhangbinlin2012@163.com}}

\subjclass[2010]{58E05, 26A33, 35J60, 47G20.}

\date{\today}
\keywords{Schr\"{o}dinger-Kirchhoff equation, Fractional magnetic operators, Variational methods}

\begin{document}

\begin{abstract}
The paper deals with existence and multiplicity of solutions
of the fractional Schr\"{o}dinger--Kirchhoff equation
involving an external magnetic potential. As a consequence,
the results can be applied to the special case
\begin{equation*}
(a+b[u]_{s,A}^{2\theta-2})(-\Delta)_A^su+V(x)u=f(x,|u|)u\,\,
\quad \text{in $\mathbb{R}^N$},
\end{equation*}
where $s\in (0,1)$, $N>2s$, $a\in \mathbb{R}^+_0$, $b\in \mathbb{R}^+_0$,
$\theta\in[1,N/(N-2s))$, $A:\mathbb{R}^N\rightarrow\mathbb{R}^N$ is a
magnetic potential, $V:\mathbb{R}^N\rightarrow \mathbb{R}^+$ is an
electric potential, $(-\Delta )_A^s$ is the fractional magnetic operator.
In the super-- and sub--linear cases, the existence of least energy solutions
for the above problem is obtained by the mountain pass theorem, combined with the
Nehari method,  and by the direct methods respectively. In the
superlinear--sublinear case, the existence of infinitely many solutions
is investigated by the symmetric mountain pass theorem.
\end{abstract}

\maketitle

\section{Introduction and main result}\label{sec1}

The paper deals with the existence of
solutions of the  fractional {\em Schr\"{o}dinger--Kirchhoff}  problem
\begin{equation}\label{eq1}
M([u]_{s,A}^2)(-\Delta)_A^su+V(x)u=f(x,|u|)u\quad \text{in $\mathbb{R}^N$},
\end{equation}
where hereafter $s\in(0,1)$, $N>2s$,
$$
[u]_{s,A}=\left(\iint_{\mathbb{R}^{2N}}\frac{|u(x)-e^{{\rm i}(x-y)\cdot A(\frac{x+y}{2})}u(y)|^2}{|x-y|^{N+2s}}dxdy\right)^{1/2},
$$
$M:\mathbb{R}^+_0\rightarrow\mathbb{R}^+_0$ is a  Kirchhoff function, $V:\mathbb{R}^N\rightarrow\mathbb{R}^+$ is a scalar
potential, $A:\mathbb{R}^N\rightarrow \mathbb{R}^N$ is a magnetic potential, and $(-\Delta )_A^s$ is the associated fractional magnetic operator which, up to a normalization constant, is defined as
\begin{equation*}
(-\Delta)_A^s\varphi(x)=2
\lim_{\varepsilon\rightarrow 0^+}\int_{\mathbb{R}^N\setminus B_\varepsilon(x)}\frac{\varphi(x)-e^{{\rm i}(x-y)\cdot A(\frac{x+y}{2})}\varphi(y)}{|x-y|^{N+2s}}\,dy,\quad x\in\mathbb{R}^{N},
\end{equation*}
along functions $\varphi\in C_0^\infty(\mathbb{R}^N,\mathbb{C})$.
Henceforward $B_\varepsilon(x)$ denotes the ball of
$\mathbb{R}^N$ centered at $x\in\mathbb{R}^N$ and
radius $\varepsilon>0$. For details on fractional magnetic operators
we refer to \cite{PDMS} and to the references \cite{ichi1,ichi2,ichi3,purice} for the physical background.

The operator $(-\Delta)_A^s$ is consistent with the definition
of fractional Laplacian $(-\Delta )^s$ when $A\equiv0$. For
further details
on $(-\Delta )^s$,  we refer the interested reader to \cite{r28}.
Nonlocal operators can be seen as the infinitesimal generators
of L\'{e}vy stable diffusion processes \cite{r8}.  Moreover, they allow
us to develop a generalization of quantum mechanics and also
to describe the motion of a chain or an array of particles
that are connected by elastic springs as well as unusual diffusion
processes in turbulent fluid motions and material transports in
fractured media (for more details see for example
\cite{r8,r3,r7} and the references therein).
Indeed, the literature on nonlocal fractional operators and
on their applications is quite large, see for example the
recent monograph \cite{MBRS}, the extensive paper \cite{DMV} and
the references cited there.

The paper was motivated by some works appeared in recent years concerning the magnetic Schr\"{o}dinger equation
\begin{align}\label{eq1.01}
-(\nabla -{\rm i} A)^2u+V(x)u=f(x,|u|)u\quad\text{in $\mathbb{R}^N,$}
\end{align}
which has been extensively studied (see \cite{GAAS,SCSS,JDJV,KK,MS}). The magnetic Schr\"{o}dinger operator is defined as
\begin{align*}
-(\nabla -{\rm i} A)^2u=-\Delta u+2{\rm i}A(x)\cdot\nabla u+|A(x)|^2u+{\rm i}u\,{\rm div}A(x).
\end{align*}
As stated in \cite{MSBV}, up to correcting the operator
by the factor $(1-s)$, it follows that
$(-\Delta)^s_A u$ converges to $-(\nabla u-{\rm i} A)^2u$
as $s\uparrow1$. Thus, up to
normalization, the nonlocal case can be seen as an approximation of the local case (see Section~\ref{singularsect} for further details).
As $A=0$ and $M=1$, equation \eqref{eq1} becomes the fractional Schr\"{o}dinger equation
\begin{align*}
(-\Delta)^su+V(x)u=f(x,|u|)u\quad\text{in $\mathbb{R}^N,$}
\end{align*}
introduced by Laskin \cite{laskin1, laskin2}. Here
the nonlinearity $f$ satisfies general conditions.
We refer, for instance, to \cite{r15, r14, r17} and the
references therein for recent results.

Throughout the paper, without explicit mention, we also assume that
{\em $A:\mathbb{R}^N\rightarrow\mathbb{R}^N$
and $V:\mathbb{R}^N\rightarrow\mathbb{R}^+$ are continuous
functions, and that $V$ satisfies},
\begin{itemize}
\item[$(V_1)$] {\em there exists $V_0>0$ such that $\inf_{\mathbb{R}^N} V \geq V_0$}.
\end{itemize}
The  Kirchhoff function $M:\mathbb{R}^+_0\rightarrow \mathbb{R}^+_0$
is assumed to be {\em continuous and to verify}
\begin{itemize}
 \item[$(M_1)$] {\em for any $\tau>0$ there exists $\kappa=\kappa(\tau)>0$ such that $M(t)\geq\kappa$ for all $t\geq\tau$};
  \item[$(M_2)$] {\em there exists $\theta\in[1,2_s^*/2)$ such that $tM(t)\leq \theta\mathscr{M}(t)$ for all $t\geq0$, $\mathscr{M}(t)=\int_0^t M(\tau)d\tau$}.
\end{itemize}
A simple typical example of $M$ is given by $M(t)=a+b\,t^{\theta-1}$
for $t\in\mathbb{R}^+_0$, where $a\in\mathbb{R}^+_0$, $b\in\mathbb{R}^+_0$
and $a+b>0$. When $M$ is of this type, problem \eqref{eq1}
is said to be {\em non--degenerate} if $a>0$, while it is
called {\em degenerate} if $a=0$.

Clearly, assumptions $(M_1)$ and $(M_2)$ cover the degenerate case.
It is worth pointing out that
the degenerate case is rather interesting and is treated in well--known
papers in Kirchhoff theory, see for example \cite{dAS}.
In the large literature on degenerate Kirchhoff problems,
the transverse oscillations of a stretched string, with nonlocal
flexural rigidity, depends continuously on the Sobolev deflection
norm of $u$ via $M(\|u\|^2)$. From a physical point of view, the fact
that $M(0)=0$ means  that the base tension of the string is zero, a
very realistic model. We refer to \cite{AFP, capu, PXZ2, XZF2, XMTZ} and the
references therein for more details in bounded domains and in the whole space.
Recent existence results of solutions for fractional non--degenerate Kirchhoff problems are given, for example, in \cite{FV, XZF1, XZG, PS}.

Assumptions $(M_1)$ and $(M_2)$ on the Kirchhoff function $M$ are enough to assure  the existence of solutions of~\eqref{eq1}.
However, to get the existence of ground states, we assume
also the further mild request
\begin{itemize}
\item[$(M_3)$] {\em there exists $m_0>0$ such that
$M(t)\geq m_0 t^{\theta-1}$ for all $t\in  [0,1]$,}
\end{itemize}
{\em where $\theta$ is the number given in $(M_2)$ when $(M_2)$
is assumed, otherwise $\theta$ is any number greater
than or equal to 1}.

Of course, $(M_3)$ is satisfied also in the model case,
even when $M(0)=0$, that is in the degenerate case.
In \cite{PXZ2}, condition $(M_3)$ was also applied to
investigate the existence of entire solutions for the stationary
Kirchhoff type equations driven by the fractional $p$--Laplacian
operator in $\mathbb{R}^N$.

Superlinear nonlinearities $f$ satisfy
\begin{itemize}
\item[$(f_1)$] {\em $f\in\mathbb{R}^N\times\mathbb{R}^+\rightarrow\mathbb{R}$ is a Carath\'{e}odory function and there exist $C>0$ and $p\in (2\theta,2_s^*)$ such that}
\begin{align*}
|f(x,t)|\leq C(1+|t|^{p-2})\quad\mbox{\em for all }
(x,t)\in\mathbb{R}^N\times\mathbb{R}^+;
\end{align*}
\item[$(f_2)$] {\em There exists a constant $\mu>2\theta$ such that
$$0<\mu F(x,t)\leq f(x,t)t^2,\quad
F(x,t)=\int_0^tf(x,\tau)\tau d\tau,$$
whenever $x\in\mathbb{R}^N$ and $t\in \mathbb{R}^+$;}
\item[$(f_3)$] {\em $f(x,t)=o(1)$ as $t\rightarrow 0^+$, uniformly for $x\in\mathbb{R}^N$};
\item[$(f_4)$] $\displaystyle{\inf_{x\in\mathbb{R}^N}F(x,1)>0}$.
\end{itemize}
A typical example of $f$, verifying $(f_1)$--$(f_4)$, is given by $f(x, |u|)=|u|^{p-2}$, with $2\theta<p<2_{s}^{*}$.
The fractional solution spaces $\Ma(\mathbb{R}^N,\mathbb{C})$
and $H_{A,V}^s(\mathbb{R}^N,\mathbb{C})$
are introduced precisely in Section~\ref{sec3}.
\vskip2pt

\noindent
We say that $u\in \Ma(\mathbb{R}^N,\mathbb{C})$ (resp.\
$u\in H_{A,V}^s(\mathbb{R}^N,\mathbb{C})$) is a (weak)
{\em solution} of \eqref{eq1}, if
\begin{align*}
\Re\bigg[M([u]_{s,A}^2)\iint_{\mathbb{R}^{2N}}\!\!&\frac{\big[u(x)-e^{{\rm i}(x-y)\cdot A(\frac{x+y}{2})}u(y)\big]\cdot\big[\overline{\varphi(x)-e^{{\rm i}(x-y)\cdot A(\frac{x+y}{2})}\varphi(y)}\big]}{|x-y|^{N+2s}} dxdy+\int_{\mathbb{R}^N}Vu\overline{\varphi} dx\bigg]\\
&\hspace{9truecm}=\Re\int_{\mathbb{R}^N}f(x,|u|)u\overline{\varphi} dx,
\end{align*}
for all $\varphi\in \Ma(\mathbb{R}^N,\mathbb{C})$
(resp.\ $\varphi\in H_{A,V}^s(\mathbb{R}^N,\mathbb{C})$).

Now we are in a position to state the first existence result.

\begin{theorem}[Superlinear case]\label{th1}
Assume that $V$ satisfies $(V_1)$, $f$ satisfies $(f_1)$--$(f_4)$
and $M$ fulfills $(M_1)$--$(M_2)$. Then  \eqref{eq1} admits
a nontrivial radial mountain pass solution
$u_0\in \Ma(\mathbb{R}^N,\mathbb{C})$.
Furthermore, if $M$ satisfies $(M_1)$--$(M_3)$, then
\eqref{eq1} has a ground state $u\in \Ma(\mathbb{R}^N,\mathbb{C})$
with positive energy.
\end{theorem}

\noindent
Sublinear nonlinearities $f$ verify
\begin{itemize}
\item[$(f_5)$] {\em There exist $q\in (1,2)$ and $a\in
L^\infty_{\rm loc}(\mathbb{R}^N)\cap L^{\frac{2}{2-q}}(\mathbb{R}^N)$
such that}
\begin{align*}
|f(x,t)|\leq a(x)t^{q-2}\quad\mbox{\em for all }
(x,t)\in\mathbb{R}^N\times\mathbb{R}^+.
\end{align*}
\item[$(f_6)$] {\em There exist $q\in (1,2)$, $\delta>0$, $a_0>0$ and
a nonempty open subset $\Omega$ of $\mathbb{R}^N$ such that}
    \begin{align*}
    |f(x,t)|\geq a_0t^{q-2}\quad\mbox{\em for all } (x,t)\in\Omega\times(0,\delta).
    \end{align*}
\end{itemize}
A typical example of $f$, verifying $(f_5)$--$(f_6)$,
is $f(x, |u|)=(1+|x|^2)^{(q-2)/2}|u|^{q-2}$ with $1<q<2$.
The second result reads as follows.

\begin{theorem}[Sublinear case]\label{th2}
Assume that $V$ satisfies $(V_1)$, $f$ satisfies $(f_5)$--$(f_6)$
and $M$ is continuous in $\mathbb R^+_0$ and satisfies
$(M_1)$ and $(M_3)$, with $\vartheta\ge1$. Then \eqref{eq1}
admits a nontrivial solution $u\in H_{A,V}^s(\mathbb{R}^N,\mathbb{C})$,
which is a ground sate of \eqref{eq1}.
\end{theorem}

\noindent
To get infinitely many solutions for equation
\eqref{eq1} in the local sublinear--superlinear case, we also assume
\begin{itemize}
\item[$(V_2)$] {\em There exists $h>0$  such that
\begin{align*}
\lim_{|y|\rightarrow\infty} {\mathscr L}^N\big(\{x\in B_h(y):V(x)\leq c\}\big)=0
\end{align*}
for all $c>0$.}
\item[$(f_7)$]
{\em $F(x,t)\geq 0$ for all $(x,t)\in\mathbb{R}^N\times\mathbb{R}^+_0$,
and there exist $q\in(1, 2)$, a nonempty open subset $\Omega$ of~$\mathbb{R}^N$
and  $a_1>0$ such that}
\begin{align*}
F(x,t)\geq a_1t^{q}\quad\mbox{\em for all }
(x,t)\in\Omega\times\mathbb{R}^+.
\end{align*}
\end{itemize}
An example of $f$, which satisfies assumptions $(f_1)$ and $(f_7)$, is
\begin{align*}
f(x,t)=(1+|x|^2)^{(q-2)/2}t^{q-2}+t^{p-2}\quad\mbox{for all }
(x,t)\in\mathbb{R}^N\times\mathbb{R}^+_0,
\end{align*}
when $1<q<2\leq 2\theta<p<2_s^*$.

\begin{theorem}[Multiplicity -- local superlinear--sublinear case]\label{th3}
 Assume that $V$ satisfies $(V_1)$--$(V_2)$, that
 $f$ fulfills  $(f_1)$ and $(f_7)$ and that $M$
is a continuous function in $\mathbb R^+_0$, verifying
$(M_1)$ and $(M_3)$, with $\vartheta\ge1$.
Then \eqref{eq1} admits a sequence $(u_k)_{k}$ of
nontrivial solutions.
\end{theorem}

\begin{remark}
{\rm $(i)$ Condition $(V_2)$, which is weaker than the
coercivity assumption: $V(x)\rightarrow \infty$ as $|x|\rightarrow \infty$,
was first proposed by {\em Bartsch} and {\em Wang} in~\cite{BW}
to overcome the lack of compactness.}

$(ii)$ {\rm To our best knowledge, Theorem~\ref{th3} is the first result
for the Schr\"{o}dinger--Kirchhoff equations involving
concave--convex nonlinearities in the fractional setting.
We also refer to \cite{XZF2} for some related multiplicity results.}
\end{remark}

The paper is organized as follows.
In Section~\ref{singularsect} we provide a few remarks about
the singular limit as $s\uparrow 1$.
In Section~\ref{sec3}, we recall some necessary definitions and
properties for the functional setting.
In Section~\ref{sec4}, we obtain some preliminary results.
In Section~\ref{sec5}, the existence of ground states of \eqref{eq1} is obtained by using the mountain pass theorem
together with the Nehari method, and by the direct
methods respectively.
In Section~\ref{sec6}, the existence of infinitely many solutions of \eqref{eq1} is obtained by using the symmetric mountain pass theorem.
\medskip

\noindent
{\bf Acknowledgements}.
Xiang Mingqi was supported by the Fundamental Research Funds for the Central Universities (No.\ 3122015L014). Patrizia Pucci and Marco Squassina are members of the {\em Gruppo Nazionale per l'Analisi Ma\-te\-ma\-ti\-ca, la Probabilit\`a e le loro Applicazioni} (GNAMPA) of the {\em Istituto Nazionale di Alta Matematica} (INdAM).
The manuscript was realized within the auspices of the INdAM -- GNAMPA Project
{\em Problemi variazionali su variet\`a Riemanniane e gruppi di Carnot} (Prot\_2016\_000421). P. Pucci was partly supported by the Italian MIUR project  {\em Variational and perturbative aspects of nonlinear differential problems} (201274FYK7).
Binlin Zhang was supported by Natural Science Foundation of Heilongjiang Province of China (No.\ A201306).

\section{Remarks on the singular limit as $s\uparrow 1$}
\label{singularsect}
\renewcommand{\S}{{\mathbb S}}

\noindent
The functional framework investigated in the paper admits a very nice consistency
property with more familiar local problems, in the singular limit
as the fractional diffusion parameter $s$
approaches $1$.
Let $\Omega$ be a nonempty open subset of $\R^N$. We denote by $L^2(\Omega,
\mathbb C)$ the Lebesgue space of complex valued functions with
summable square, endowed with the norm $\|u\|_{L^2(\Omega,
\mathbb C)}$.\ We indicate by $H^s_A(\Omega)$ the space of functions
$u\in L^2(\Omega,\mathbb C)$ with finite magnetic Gagliardo semi--norm,
given by
$$
[u]_{H^s_A(\Omega)}=\left(\iint_{\Omega\times\Omega}\frac{|u(x)-e^{{\rm i}(x-y)\cdot A(\frac{x+y}{2})}u(y)|^2}{|x-y|^{N+2s}}dxdy\right)^{1/2}.
$$
The space $H^s_A(\Omega)$ is equipped with the norm
$$
\|u\|_{H^s_A(\Omega)}=
\big(\|u\|_{L^2(\Omega,\mathbb C)}^2+[u]_{H^s_A(\Omega)}^2\big)^{1/2}.
$$
The space $H^s_{0,A}(\Omega)$ is the completion of $C^\infty_c(\Omega,
\mathbb C)$ in $H^s_A(\Omega)$.

Indeed, in the recent paper~\cite{MSBV}, the following theorem was proved, which is a Bourgain--Brezis--Mironescu
type result in the framework of  magnetic Sobolev spaces.
\begin{proposition}[Theorems 1.1 and 1.2 of~\cite{MSBV}]\label{main}
Let $\Omega$ be an open bounded subset of $\R^N$, with Lipschitz
boundary and let $A$ be of class $C^2$ over $\overline{\Omega}$.
Then,
$$\lim_{s\uparrow 1}(1-s)
\iint_{\Omega\times\Omega}\frac{|u(x)-e^{\i (x-y)
\cdot A\left(\frac{x+y}{2}\right)}u(y)|^2}{|x-y|^{N+2s}}dxdy
=K_N\int_{\Omega}|\nabla u-\i A(x)u|^2dx
$$
for every  $u\in H^1_{A}(\Omega)$, where
\begin{align*}
	K_{N}=\frac{1}{2}\int_{\S^{N-1}}|\omega\cdot
{\bf e}|^{2}d\mathcal{H}^{N-1}(\omega),
\end{align*}
and $\S^{N-1}$ is the unit sphere of $\R^N$
and ${\bf e}$ any unit vector of $\R^{N}$. Furthermore,
$$\lim_{s\uparrow 1}(1-s)
\iint_{\R^{2N}}\frac{|u(x)-e^{\i (x-y)\cdot
A\left(\frac{x+y}{2}\right)}u(y)|^2}{|x-y|^{N+2s}}dxdy=
K_N\int_{\Omega}|\nabla u-\i A(x)u|^2dx
$$
for every  $u\in H^1_{0,A}(\Omega)$.
\end{proposition}

Problem \eqref{eq1} could be treated in an arbitrary smooth
open bounded subset $\Omega$ of~$\R^N$, provided that
the solution space
is~$W$, which consists of all functions $u$ in $H^s_{A}(\R^N)$,
with $u=0$ in $\R^N\setminus\Omega$.  More precisely, consider
the non--degenerate model case
$$
M(t)=a(s)+b(s)t,\qquad\mbox{where
$a(s)\approx 1-s\,\,$ and\,\, $b(s)\approx (1-s)^2b_0$\quad as $s\uparrow 1$.}
$$
Then the corresponding problem \eqref{eq1} in $\Omega$ writes as
\begin{equation*}
\begin{cases}
\left(1+(1-s)b_0\displaystyle\iint_{\R^{2N}}
\frac{|u(x)-e^{\i (x-y)\cdot A\left(\frac{x+y}{2}\right)}u(y)|^2}
{|x-y|^{N+2s}}dxdy\right)\widehat{(-\Delta)_A^s} u+V(x)u=f(x,|u|)u&\mbox{in }\Omega,\\
u=0  &\mbox{in }\mathbb{R}^N\setminus\Omega,
\end{cases}
\end{equation*}
where $u$ belongs to the solution space  $W$ and
$$
\widehat{(-\Delta)_A^s} u=(1-s)(-\Delta)_A^su.
$$
This is natural since the Gagliardo semi--norms are typically multiplied by normalizing constants which vanish
at the rate of $1-s.$ Since by Proposition~\ref{main}
\begin{align*}
(1-s)\iint_{\R^{2N}}\frac{|u(x)-e^{\i (x-y)\cdot A\left(\frac{x+y}{2}\right)}u(y)|^2}{|x-y|^{N+2s}}dxdy & \approx
\int_{\Omega}|\nabla u-\i A(x)u|^2dx\quad\text{as }s\uparrow 1,\\
\widehat{(-\Delta)_A^s} u =(1-s)(-\Delta)_A^su &\approx -(\nabla u-{\rm i} A)^2u,\quad\text{as }s\uparrow 1,
\end{align*}
the above problem converges to the local problem
\begin{equation*}
\begin{cases}
-\left(1+ b_0\displaystyle\int_{\Omega}|\nabla u-\i A(x)u|^2dx\right)(\nabla u-{\rm i} A)^2u+V(x)u=f(x,|u|)u &  \mbox{in $\Omega$},\\
u=0  &  \mbox{on $\partial\Omega$},
\end{cases}
\end{equation*}
which as $A\to O$ reduces  to
\begin{equation*}
\begin{cases}
-\left(1+ b_0\displaystyle\int_{\Omega}|\nabla u|^2dx\right)
\Delta u+V(x)u=f(x,|u|)u &  \mbox{in }\Omega,\\
u=0  &  \mbox{on }\partial\Omega.
\end{cases}
\end{equation*}
This is the classical model of a Schr\"{o}dinger--Kirchhoff
equation. When $b_0=0$, the last two problems become the
classical Schr\"odinger Dirichlet problems with or
without external magnetic potential~$A$.

\section{Functional setup}\label{sec3}
We first provide some basic functional setting that will be used in the next sections.
The critical exponent $2^*_s$ is defined as $2N/(N-2s).$
Let $L^2(\mathbb{R}^N,V)$ denote the Lebesgue space of real valued functions with
$V(x)|u|^2\in L^1({\mathbb{R}}^N),$ equipped with norm
$$
\|u\|_{2,V}=\left(\int_{\mathbb{R}^N}V(x)|u|^2 dx\right)^{{1}/{2}}\quad \text{for all }u\in L^2(\mathbb{R}^N,V).
$$
The fractional Sobolev space $H_{V}^{s}(\mathbb{R}^N)$ is then defined as
\begin{align*}
H^{s}_V(\mathbb{R}^N)=\big\{u\in L^2(\mathbb{R}^N,V):[u]_{s}<\infty\big\},
\end{align*}
where $[u]_{s}$ is the Gagliardo semi--norm
\begin{align*}
[u]_{s}=\left(\iint_{\mathbb{R}^{2N}}\frac{|u(x)-u(y)|^2}
{|x-y|^{N+2s}}dxdy\right)^{{1}/{2}}.
\end{align*}
The space $H_{V}^{s}(\mathbb{R}^N)$ is endowed  with the norm
\begin{align*}
\|u\|_s=\left(\|u\|_{2,V}^2+[u]_{s}^2\right)^{{1}/{2}}.
\end{align*}
The localized norm, on a compact subset $K$ of $\mathbb R^N$,
for the space $H^s_V(K)$, is denoted by
\begin{equation}
\label{localiz}
\|u\|_{s,K}=\left(\int_K V(x)|u|^2dx+\iint_{K\times
K}\frac{|u(x)-u(y)|^2}{|x-y|^{N+2s}}dxdy\right)^{{1}/{2}}.
\end{equation}
The embedding $H^{s}_V(\mathbb{R}^N)\hookrightarrow L^{\nu}(\mathbb{R}^N)$ is continuous for any $\nu\in [2,2_s^*]$ by
\cite[Theorem 6.7]{r28}, namely there exists a positive constant $C$ such that
\begin{align*}
\|u\|_{L^{\nu}(\mathbb{R}^N)}\leq C \|u\|_s\quad\mbox{for all } u\in H^{s}_V(\mathbb{R}^N).
\end{align*}
Let us set
$$
H_{r,V}^{s}(\mathbb{R}^N)=\big\{u\in H^{s}_V(\mathbb{R}^N):u(x)=u(|x|)
\mbox{ for all }x\in\mathbb{R}^N\big\}.
$$
To prove the existence of radial weak solutions of \eqref{eq1}, we shall
use the following embedding theorem due to P.L. Lions.

\begin{theorem}[Compact embedding, I -- Th\'{e}or\`{e}me II.1
of \cite{Lions}]\label{th2.1} Let $N\ge2$.
For any $\alpha\in(2,2_s^*)$ the  embedding
$
H_{r,V}^{s}(\mathbb{R}^N)\hookrightarrow\hookrightarrow L^{\alpha}(\mathbb{R}^N)
$
is compact.
\end{theorem}

\noindent
Furthermore, we also have

\begin{theorem}[Compact embedding, II -- Theorem 2.1 of \cite{PXZ}]
\label{th2.2}
Assume that conditions $(V_1)$--$(V_2)$ hold. Then, for any $\nu\in(2,2_s^*)$
the embedding
$
H_V^s(\mathbb{R}^N)\hookrightarrow\hookrightarrow  L^\nu(\mathbb{R}^N)
$
is compact.
\end{theorem}

\noindent
Let $L_V^2(\mathbb{R}^N,\mathbb{C})$ be the Lebesgue space of functions $u:\mathbb{R}^N\to\mathbb{C}$
with $V|u|^2 \in L^1({\mathbb R}^N)$, endowed with the (real) scalar product
\begin{align*}
\langle u,v\rangle_{L^2,V}=\Re \int_{\mathbb{R}^N}V(x)u \overline{v}dx\quad \mbox{for all } u,v\in L^2(\mathbb{R}^N,\mathbb{C}),
\end{align*}
where $\bar z$ denotes complex conjugation of $z\in \mathbb{C}$.
Consider now, according to \cite{PDMS}, the magnetic Gagliardo
semi--norm given by
\begin{align*}
[u]_{s,A}=\left(\iint_{\mathbb{R}^{2N}}\frac{|u(x)-e^{{\rm i}(x-y)
\cdot A(\frac{x+y}{2})}u(y)|^2}{|x-y|^{N+2s}}dxdy\right)^{1/2}.
\end{align*}
Define $H_{A,V}^{s}(\mathbb{R}^N)$ as the closure of $C_c^\infty(\mathbb{R}^N,\mathbb{C})$ with respect to the norm
\begin{align*}
\|u\|_{s,A}=\big(\|u\|_{L^2,V}^2+[u]_{s,A}^2\big)^{1/2}.
\end{align*}
A scalar product on $H_{A,V}^{s}(\mathbb{R}^N)$ is given by
\begin{align*}
\langle u,v\rangle_{s,A}=\langle u,v\rangle_{L^2,V}+\Re \iint_{\mathbb{R}^{2N}}\frac{\big[u(x)-e^{{\rm i}(x-y)\cdot A(\frac{x+y}{2})}u(y)\big]\cdot\big[\overline{v(x)-e^{{\rm i}(x-y)\cdot A(\frac{x+y}{2})}v(y)}\big]}{|x-y|^{N+2s}} dxdy.
\end{align*}
Arguing as in \cite[Proposition 2.1]{PDMS}, we see
that $\big(H_{A,V}^{s}(\mathbb{R}^N),\langle\cdot,\cdot\rangle_{s,A}\big)$
is a real Hilbert space.

\begin{lemma}\label{lemma2.1}
For each $u\in H^s_{A,V}(\mathbb{R}^N,\mathbb{C})$
\begin{align*}
|u|\in H^s_V(\mathbb{R}^N)\quad\mbox{and}\quad
\big\||u|\big\|_{s}\leq \|u\|_{s,A}.
\end{align*}
\end{lemma}

\begin{proof}
The assertion follows directly from the pointwise diamagnetic inequality
\begin{align*}
\big||u(x)|-|u(y)|\big|\leq \left|u(x)-e^{{\rm i}(x-y)\cdot A(\frac{x+y}{2})}u(y)\right|,
\end{align*}
 for a.e. $x,y\in\mathbb{R}^N$, see \cite[Lemma 3.1, Remark 3.2]{PDMS}.
\end{proof}

\noindent
Following Lemma~\ref{lemma2.1} and  using the same discussion of~\cite[Lemma 3.5]{PDMS}, we have

\begin{lemma}\label{lemma2.2}
The embedding
$$
H_{A,V}^s(\mathbb{R}^N,\mathbb{C})\hookrightarrow L^p(\mathbb{R}^N,\mathbb{C})
$$
is continuous for all $p\in[2,2_s^*]$. Furthermore, for any compact
subset $K\subset \mathbb{R}^N$  and all $p\in[1,2_s^*)$ the embeddings
$$
H_{A,V}^s(\mathbb{R}^N,\mathbb{C})\hookrightarrow
H^s_V(K,\mathbb{C})\hookrightarrow \hookrightarrow  L^p(K,\mathbb{C})
$$
are continuous and the latter is compact, where
$H^s_V(K,\mathbb{C})$ is endowed with \eqref{localiz}.
\end{lemma}

\noindent
Define now
\begin{align*}
\Ma(\mathbb{R}^N,\mathbb{C})=\big\{u\in H_{A,V}^s(\mathbb{R}^N,\mathbb{C}):u(x)=u(|x|),\,\, x\in\mathbb{R}^N\big\}.
\end{align*}
By Theorems~\ref{th2.1}--\ref{th2.2}  and Lemma~\ref{lemma2.1}, we have the following lemma (cf. also \cite[Lemma 4.1]{PDMS}).

\begin{lemma}\label{lemma2.3}
Let $V$ satisfy $(V_1)$. Let $(u_n)_{n}$ be a bounded sequence
in $\Ma(\mathbb{R}^N,\mathbb{C})$. Then, up to a subsequence,
$(|u_n|)_{n}$ converges strongly to some function $u$
in $L^p(\mathbb{R}^N)$ for all $p\in(2,2_s^*)$.

Moreover, if $V$ satisfies $(V_1)$--$(V_2)$, then for all bounded
sequence $(u_n)_{n}$ in $H_{A,V}^s(\mathbb{R}^N,\mathbb{C})$
the sequence $(|u_n|)_{n}$ admits a subsequence converging strongly
to some $u$ in $L^p(\mathbb{R}^N)$ for all $p\in[2,2_s^*)$.
\end{lemma}

\section{Preliminary results}\label{sec4}

The functional $\mathcal{I}:\Ma(\mathbb{R}^N,\mathbb{C})\to\mathbb{R}$, associated
with equation \eqref{eq1}, is defined by
\begin{align*}
\mathcal{I}(u)=\frac{1}{2}\mathscr{M}([u]_{s,A}^2)
+\frac{1}{2}\|u\|_{L^2,V}^2-\int_{\mathbb{R}^N}F(x,|u|)dx.
\end{align*}
It is easy to see that $\mathcal{I}$ is of class $C^1(\Ma(\mathbb{R}^N,\mathbb{C}),\mathbb{R})$ and
\begin{align*}
\langle \mathcal{I}^\prime(u),v\rangle=
&\Re\bigg[M([u]_{s,A}^2)\iint_{\mathbb{R}^{2N}}\frac{(u(x)-e^{{\rm i}(x-y)\cdot A(\frac{x+y}{2})}u(y))(\overline{v(x)-e^{{\rm i}(x-y)\cdot A(\frac{x+y}{2})}v(y)})}{|x-y|^{N+2s}} dxdy\\
&\hspace{6truecm}+\int_{\mathbb{R}^N}Vu\overline{v} dx\bigg]-\Re\int_{\mathbb{R}^N}f(x,|u|)u\overline{v}dx,
\end{align*}
for all $u,v\in \Ma(\mathbb{R}^N,\mathbb{C})$. Hereafter, $\langle\cdot,\cdot\rangle$ denotes the
duality pairing between $\big(\Ma(\mathbb{R}^N,\mathbb{C})\big)'$ and
$\Ma(\mathbb{R}^N,\mathbb{C})$.

Hence, the critical points of $\mathcal{I}$ are exactly the weak solutions of
\eqref{eq1}. Moreover, $\mathscr{M}([u]_{s,A}^2)$
is weakly lower semi--continuous in $\Ma(\mathbb{R}^N,\mathbb{C})$
 by the weak lower semi--continuity of
$u\mapsto [u]_{s,A}^2$ jointly with the monotonicity
and continuity of $\mathscr{M}$. Hence, $\mathcal{I}$
is weakly lower semi--continuous in $\Ma(\mathbb{R}^N,\mathbb{C})$,
being $\int_{\mathbb{R}^N}F(x,|u|)dx$ weakly continuous
in $\Ma(\mathbb{R}^N,\mathbb{C})$.

\begin{definition} {\rm We say that $\mathcal{I}$ satisfies the $(PS)$ {\em condition} in $\Ma(\mathbb{R}^N,\mathbb{C})$, if any $(PS)$
	sequence $(u_n)_{n}\subset \Ma(\mathbb{R}^N,\mathbb{C})$, namely
a sequence such that $(\mathcal{I}(u_n))_{n}$ is bounded and
$\mathcal{I}^\prime (u_n)\rightarrow 0$ as $n\rightarrow\infty$, admits a strongly convergent
subsequence in $\Ma(\mathbb{R}^N,\mathbb{C})$.}
\end{definition}

\begin{lemma}[Palais--Smale condition]
	\label{lemma3.1}
Let $(M_1)$--$(M_2)$ and $(f_1)$--$(f_3)$ hold. Then $\mathcal{I}$
satisfies the $(PS)$ condition in
$\Ma(\mathbb{R}^N,\mathbb{C})$.
\end{lemma}

\begin{proof}
Let $(u_n)_{n}$ be a $(PS)$ sequence in $\Ma(\mathbb{R}^N,\mathbb{C})$. Then there
exists $C>0$ such that
$|\mathcal{I}(u_n)|\leq C$ and $|\langle \mathcal{I}^\prime (u_n),u_n\rangle|\leq C\|u_n\|_{s,A}$
for all $n$. As in Lemma~4.5 of~\cite{capu}, see
also \cite{colpuc},
we divide the proof into two parts.

\smallskip
\noindent
$\bullet$\,{\em Case} $\inf_{n\in\mathbb N}[u_n]_{s,A}=d>0$.
By $(M_1)$, there
exists $\kappa=\kappa(d)>0$ with $M(t)\geq \kappa>0$ for all $t\geq d$.
Thus, $(M_2)$ and $(f_2)$ yield
\begin{equation}\label{NEW}\begin{aligned}
C+C\|u_n\|_{s,A}&\geq \mathcal{I}(u_n)-\frac{1}{\mu}
\langle \mathcal{I}^\prime (u_n),u_n\rangle\\
&=\frac{1}{2}\mathscr{M}([u_n]_{s,A}^2)-\frac{1}{\mu}M([u_n]_{s,A}^2)[u_n]_{s,A}^2
+\left(\frac{1}{2}-\frac{1}{\mu}\right)\|u_n\|_{L^2,V}^2\\
&\ \ -\frac{1}{\mu}\int_{\mathbb{R}^N}(\mu F(x,|u_n|)-f(x,|u_n|)|u_n|^2)dx\\
&\geq \frac{1}{2}\mathscr{M}([u_n]_{s,A}^2)-\frac{1}{\mu}M([u_n]_{s,A}^2)[u_n]_{s,A}^2
+\left(\frac{1}{2}-\frac{1}{\mu}\right)\|u_n\|_{L^2,V}^2\\
&\geq \left(\frac{1}{2\theta}-\frac{1}{\mu}\right)M([u_n]_{s,A}^2)[u_n]_{s,A}^2
+\left(\frac{1}{2}-\frac{1}{\mu}\right)\|u_n\|_{L^2,V}^2\\
&\geq \kappa\left(\frac{1}{2\theta}-\frac{1}{\mu}\right)[u_n]_{s,A}^2+\left(\frac{1}{2}
-\frac{1}{\mu}\right)\|u_n\|_{L^2,V}^2.
\end{aligned}\end{equation}
This implies at once that $(u_n)_{n}$ is bounded in $\Ma(\mathbb{R}^N,\mathbb{C})$, being $\mu>2\theta$. Going if necessary to a subsequence,
thanks to Lemmas~\ref{lemma2.2} and~\ref{lemma2.3}, we have
\begin{align}\label{eq3.1}
&u_n\rightharpoonup u\ \ {\rm in}\ \Ma(\mathbb{R}^N,\mathbb{C}),\ \ u_n\rightarrow u\ \ {\rm a.e.\ in}\ \mathbb{R}^N,\nonumber\\
&|u_n|\rightarrow |u|\ \ {\rm in}\ L^p(\mathbb{R}^N), \\
&|u_n|\leq h\ \ {\rm a.e.\ in}\ \mathbb{R}^N,\ \ {\rm for\ some}\
h\in L^p(\mathbb{R}^N). \nonumber
\end{align}
To prove that $(u_n)_{n}$ converges strongly  to $u$ in $\Ma(\mathbb{R}^N,\mathbb{C})$ as $n\to\infty$,
 we first introduce a simple notation. Let $\varphi\in \Ma(\mathbb{R}^N,\mathbb{C})$ be fixed and
denote by $L(\varphi)$ the linear functional on $\Ma(\mathbb{R}^N,\mathbb{C})$ defined by
\begin{align}\label{Lu}
\langle L(\varphi), v\rangle
=\Re\iint_{\mathbb{R}^{2N}}
\frac{(\varphi(x)-e^{{\rm i}(x-y)\cdot
 A(\frac{x+y}{2})}\varphi(y))}
 {|x-y|^{N+2s}}(\overline{v(x)-e^{{\rm i}(x-y)\cdot
 A(\frac{x+y}{2})}v(y)})dxdy,
\end{align}
for all $v\in\Ma(\mathbb{R}^N,\mathbb{C})$.  Clearly, by the H\"{o}lder inequality,  $L(\varphi)$ is continuous, being
\begin{align*}
|\langle L(\varphi),v\rangle|
\leq\|\varphi\|_{s,A}\|v\|_{s,A}.
\end{align*}
Hence the weak convergence in \eqref{eq3.1} gives
\begin{align*}
\lim_{n\rightarrow\infty}\langle L(u),u_n-u\rangle=0.
\end{align*}
Further, by the boundedness of $M([u_n]_{s,A}^2)$ we have
\begin{align}\label{eq3.2}
\lim_{n\rightarrow\infty}M([u_n]_{s,A}^2)\langle L(u),u_n-u\rangle=0.
\end{align}
By $(f_1)$ and $(f_3)$, for any $\varepsilon>0$ there exists $C_\varepsilon>0$ such that
\begin{equation}\label{g2}
|f(x,t)t|\leq \varepsilon|t|+C_\varepsilon|t|^{p-1}\ \ {\rm for \ all}\ x\in\mathbb{R}^N
\ {\rm and}\ t\in\mathbb{R}^+.
\end{equation}
Using the H\"{o}lder inequality, we obtain
\begin{equation}\label{eq3.3}\begin{aligned}
\int_{\mathbb{R}^N}\big|(f(x,&|u_n|)u_n-f(x,|u|)u)(\overline{u_n-u})\big|dx\\
&\leq \int_{\mathbb{R}^N}
[\varepsilon(|u_n|+|u|)+C_\varepsilon(|u_n|^{p-1}+|u|^{p-1})]|u_n-u|dx\\
&\leq \varepsilon(\|u_n\|_{L^2}+\|u\|_{L^2})\|u_n-u\|_{L^2}
+C_\varepsilon(\|u_n\|_{L^p}^{p-1}+\|u\|_{L^p}^{p-1})\|u_n-u\|_{L^p}\\
&\leq C\varepsilon+CC_\varepsilon\|u_n-u\|_{L^p}.
\end{aligned}\end{equation}
The Brezis--Lieb lemma and the fact that $|u_n|\rightarrow |u|$ in $L^p(\mathbb{R}^N)$ give
\begin{align*}
\lim_{n\rightarrow\infty}\int_{\mathbb{R}^N}|u_n-u|^pdx
=\lim_{n\rightarrow\infty}\int_{\mathbb{R}^N}\big(|u_n|^p-|u|^p
\big)dx=0.\end{align*}
Inserting this in \eqref{eq3.3}, we get
\begin{align}\label{eq3.4}
\lim_{n\rightarrow\infty}
\int_{\mathbb{R}^N}(f(x,|u_n|)u_n-f(x,|u|)u)(\overline{u_n-u})dx=0,
\end{align}
since $\varepsilon$ is arbitrary.
Of course, $\langle\mathcal{I}'(u_n)-\mathcal{I}'(u),u_n-u\rangle\rightarrow0$ as $n\to\infty$, since $u_n\rightharpoonup u$ in $\Ma(\mathbb{R}^N,\mathbb{C})$ and $\mathcal{I}^\prime(u_n)\to 0$ in the dual space of $\Ma(\mathbb{R}^N,\mathbb{C})$. Thus,
\begin{align*}
o(1)&=\langle\mathcal{I}^\prime(u_n)-\mathcal{I}^\prime(u),u_n-u\rangle\\
&=M([u_n]_{s,A}^2)\langle L(u_n)-L(u),u_n-u\rangle+\|u_n-u\|_{L^2,V}^2\\
&+\left(M([u_n]_{s,A}^2)-M([u]_{s,A}^2)\right)\langle L(u),u_n-u\rangle-\Re\int_{\mathbb{R}^N}(f(x,|u_n|)u_n-f(x,|u|)u)(\overline{u_n-u})dx,
\end{align*}
this, together with \eqref{eq3.2} and \eqref{eq3.4}, implies that
\begin{align*}\lim_{n\rightarrow\infty}\left(
M([u_n]_{s,A}^2)\langle L(u_n)-L(u),u_n-u\rangle+\|u_n-u\|_{L^2,V}^2
\right)=0,
\end{align*}
 which yields $u_n\rightarrow u$ in
 $\Ma(\mathbb{R}^N,\mathbb{C})$, since $M([u_n]_{s,A}^2)\ge\kappa>0$
 for all $n\geq 1$.
\smallskip

\noindent
$\bullet$\,\textit{Case} $\inf_{n\in\mathbb N}[u_n]_{s,A}=0$.
If $0$ is an isolated point for $([u_n]_{s,A})_n$, then there is a subsequence $([u_{n_k}]_{s,A})_k$ such that
$\inf_{k\in\mathbb N}[u_{n_k}]_{s,A}=d>0$ and one can proceed as before.\
If, instead, $0$ is an accumulation point for $([u_n]_{s,A})_n$, there is a subsequence, still labeled as $(u_n)_n$, such that
\begin{equation}\label{con5}
[u_n]_{s,A}\to 0,\,\,\,
\mbox{$u_n\to 0$ in }L^{2^*_s}(\mathbb R^N)\mbox{ and a.e. in }\mathbb R^N.
\end{equation}
We claim that $(u_n)_n$ converges strongly to $0$ in
$\Ma(\mathbb{R}^N,\mathbb{C})$. To this aim, we need only to show that
$\|u_n\|_{2,V}\to 0$ thanks to \eqref{con5}.
Now, \eqref{NEW} and \eqref{con5} yield that as $n\to\infty$
\begin{equation*}
C+C\|u_n\|_{2,V}+o(1)\ge \left(\frac{1}{2}
-\frac{1}{\mu}\right)\|u_n\|^2_{2,V}+o(1).
\end{equation*}
Hence, $(u_n)_n$ is bounded in $L^2(\mathbb R^N,V)$ and so
in $\Ma(\mathbb{R}^N,\mathbb{C})$. Thus, by \eqref{con5}
and  Lemma~\ref{lemma2.2}
\begin{align}\label{con6}
u_n\rightharpoonup 0\mbox{ in $\Ma(\mathbb{R}^N,\mathbb{C})$
and $u_n\rightarrow 0$ in }L^{p}(\mathbb R^N),\end{align}
being $p\in(2,2^*_s)$.
Clearly, by \eqref{g2} and \eqref{con6}, for every $\varepsilon>0$
$$
\Bigg|\int_{\mathbb R^N}f(x,|u_n|)u_n^2dx\Bigg|\le\varepsilon\|u_n\|^{2}_{2} +C_\varepsilon\|u_n\|^p_p=\varepsilon C+o(1)
$$
as $n\to\infty$. Thus,
\begin{equation}\label{x9}
\lim_{n\to\infty}\int_{\mathbb R^N}f(x,|u_n|)u_n^2dx=0,
\end{equation}
being $\varepsilon>0$ arbitrary.
Obviously, $\langle\mathcal I'(u_n),u_n\rangle\to 0$ as $n\to\infty$, by \eqref{con6} and the fact that $\mathcal I'(u_n)\to 0$ in $\big(\Ma(\mathbb{R}^N,\mathbb{C})\big)'$. Hence, by the continuity of $M$ and \eqref{con5}--\eqref{x9}, we have
\begin{align*}
o(1)&=\langle \mathcal
I'(u_n),u_n\rangle=M([u_n]_{s,A}^2)[u_n]_{s,A}^2+\|u_n\|_{2,V}^2
-\int_{\mathbb R^N}f(x,|u_n|)u_n^2dx\\
&=\|u_n\|_{2,V}^2+o(1)
\end{align*}
as $n\to\infty$. This shows the claim.

Therefore, $\mathcal I$ satisfies the $(PS)$ condition
in $\Ma(\mathbb{R}^N,\mathbb{C})$ also in this second
case and this completes the proof.
\end{proof}

\noindent
Before going to the proof of Theorem \ref{th1}, we give some useful preliminary results.

\begin{lemma}[Mountain Pass Geometry I]\label{lemma3.2}
Assume that $(M_1)$--$(M_2)$, $(f_1)$ and $(f_3)$ hold. Then there exist constant $\varrho,\alpha>0$ such that
$\mathcal{I}(u)\geq \alpha$ for all $u\in \Ma(\mathbb{R}^N,\mathbb{C})$ with $\|u\|_{s,A}=\varrho$.
\end{lemma}

\begin{proof}
It follows from $(f_3)$ that for any $\varepsilon\in (0,1)$ there
exists $\delta=\delta(\varepsilon)>0$ such that
$|f(x,t)|\leq \varepsilon$  for all $x\in\mathbb{R}^N$
and $t\in[0,\delta]$.
On the other hand, $(f_1)$ yields that
$|f(x,t)|\leq C\big(1+\delta^{2-p}\big)|t|^{p-2}$
for all $x\in\mathbb{R}^N$ and $t>\delta$. In conclusion,
\begin{align}\label{eqbu}
|f(x,t)|\leq \varepsilon+C\big(1+\delta^{2-p}\big)|t|^{p-2}\quad
\mbox{for  all }  x\in\mathbb{R}^N \mbox{ and } t\in\mathbb R^+_0.
\end{align}
Whence, for some $C_\eps>0$, we get
\begin{align}\label{eq3.5}
|F(x,t)|\leq \int_0^t |f(x,\tau)\tau| d\tau
\leq \frac{\varepsilon}{2}t^2+C_\varepsilon t^{p},
\end{align}
for all $x\in\mathbb{R}^N$  and  $t\geq0$.
Moreover, $(M_2)$ gives
\begin{align}\label{eqbu3.7}
\mathscr{M}(t)\geq \mathscr{M}(1)t^\theta\quad \mbox{for all }
t\in[0, 1],
\end{align}
while $(M_1)$ implies that $ \mathscr{M}(1)>0$.
Thus, using \eqref{eq3.5}, \eqref{eqbu3.7} and the H\"{o}lder
inequality, we obtain
for all $u\in \Ma(\mathbb{R}^N,\mathbb{C})$, with $\|u\|_{s,A}\leq 1$,
\begin{align*}
\mathcal{I}(u)&=\frac{1}{2}\mathscr{M}(\|u\|_{s,A}^2)+\frac{1}{2}\|u\|_{L^2,V}^2
-\int_{\mathbb{R}^N}F(x,|u|)dx\\
&\geq \frac{\mathscr{M}(1)}{2}[u]_{s,A}^{2\theta}+\frac{1}{2}\|u\|_{L^2,V}^2
-\frac{\varepsilon}{2}\int_{\mathbb{R}^N}|u|^2dx
-C_\varepsilon\int_{\mathbb{R^N}}|u|^pdx\\
&\geq \frac{\mathscr{M}(1)}{2}[u]_{s,A}^{2\theta}
+\Big(\frac{1}{2}-\frac{\varepsilon}{2V_0}\Big)\|u\|_{L^2,V}^2
-C_\varepsilon C_p^p\|u\|_{s,A}^p\\
&\geq \min\left\{\frac{\mathscr{M}(1)}{2},\frac{V_0-\varepsilon}{2V_0}\right\}([u]_{s,A}^{2\theta}
+\|u\|_{L^2,V}^2)-C_\varepsilon C_p^p\|u\|_{s,A}^p\\
&\geq \min\left\{\frac{\mathscr{M}(1)}{2},\frac{V_0-\varepsilon}{2V_0}\right\}([u]_{s,A}^{2\theta}
+\|u\|_{L^2,V}^{2\theta})-C_\varepsilon C_p^p\|u\|_{s,A}^p\\
&\geq 2^{1-\theta}\min\left\{\frac{\mathscr{M}(1)}{2},\frac{V_0-\varepsilon}{2V_0}\right\}([u]_{s,A}^{2}
+\|u\|_{L^2,V}^2)^{\theta}-C_\varepsilon C_p^p\|u\|_{s,A}^p\\
&=\left(2^{1-\theta}\min\left\{\frac{\mathscr{M}(1)}{2},
\frac{V_0-\varepsilon}{2V_0}\right\}-
C_\varepsilon C_p^p\|u\|_{s,A}^{p-2\vartheta}\right)
\|u\|_{s,A}^{2\vartheta},
\end{align*}
where $C_p$ is the embedding constant of $\Ma(\mathbb{R}^N,\mathbb{C})$ into
$L^p(\mathbb{R}^N,\mathbb{C})$ given by Lemma~\ref{lemma2.2}. Here we used that $\|u_n\|_{L^2,V}\leq \|u_n\|_{s,A}\leq 1$ and the inequality $(a+b)^\theta\leq 2^{\theta-1}(a^\theta+b^\theta)$ for all $a,b\geq0$.
Choosing $\varepsilon=V_0/2$ and taking $\|u\|_{s,A}=\varrho\in(0,1)$
so small that
$$
2^{1-\theta}\min\left\{\frac{\mathscr{M}(1)}{2},\frac{V_0}{4}\right\}-C_{V_0/2} C_p^p\varrho^{p-2\vartheta}>0,
$$
we have
\begin{align*}
\mathcal{I}(u)\geq \alpha=
\left(2^{1-\theta}\min\left\{\frac{\mathscr{M}(1)}{2},
\frac{V_0}{4}\right\}-C_{V_0/2} C_p^p\varrho^{p-2\vartheta}\right)
\varrho^{2\vartheta}>0,
\end{align*}
for all $u\in \Ma(\mathbb{R}^N,\mathbb{C})$, with $\|u\|_{s,A}=\varrho$.
\end{proof}

\begin{lemma}[Mountain Pass Geometry II]\label{lemma3.3}
Assume that $(M_1)$--$(M_2)$ and $(f_1)$--$(f_4)$ hold.
Then there exists
$e\in C_c^\infty(\mathbb{R}^N,\mathbb{C})$, with
$\|e\|_{s,A}\ge2$, such that
$\mathcal{I}(e)<0$. In particular, $\|e\|_{s,A}>\rho$,
where $\rho>0$ is the number introduced in Lemma~$\ref{lemma3.2}$
\end{lemma}
\begin{proof}

For any $x\in \mathbb{R}^N$, set $k(t)=F(x,t)t^{-\mu}$ for all $t\geq 1$.
Condition $(f_2)$ implies
that $k$ is nondecreasing on $[1,\infty)$. Therefore,
$k(t)\geq k(1)$ for any $t\geq1$,  that is,
\begin{align}\label{eq3.6}
F(x,t)\geq F(x,1)t^\mu\geq c_F|t|^\mu\quad \mbox{for all }x\in\mathbb{R}^N
\mbox{and }t\geq1,
\end{align}
where $c_F=\inf_{x\in\mathbb{R}^N}F(x,1)>0$ by assumption $(f_4)$.
{From} $(f_3)$ there exists $\delta\in(0,1)$ such that
$|f(x,t)t|\leq t$ for all $x \in\mathbb{R}^N$
and  $t\in[0,\delta]$.
Furthermore,  $|f(x,t)|\leq 2C$ for all $x\in\mathbb{R}^N$ and all $t$, with $\delta<t\leq1$,  thanks to $(f_1)$.
Hence, the above inequalities imply that
$f(x,t)t\geq -(1+2C)t$ for $x\in\mathbb{R}^N$ and $t\in[0,1]$.
Thus,
\begin{align}\label{eq3.7}
F(x,t)=\int_0^t f(x,\tau)\tau d\tau\geq -\frac{1+2C}{2}t^2\quad \mbox{for
all } x\in\mathbb{R}^N\mbox{ and }t\in[0,1].
\end{align}
Combining \eqref{eq3.6} with \eqref{eq3.7}, we obtain
\begin{align}\label{eq3.8}
F(x,t)\geq c_F|t|^\mu-C_F|t|^2\ \ {\rm for\ all}\ x\in\mathbb{R}^N\ {\rm and}\ t\geq0,
\end{align}
where $C_F=c_F+(1+2C)/2$.
Again $(M_2)$ gives
\begin{align}\label{eqbu3.8}
\mathscr{M}(t)\leq \mathscr{M}(1)t^\theta\quad \mbox{for all }t\geq1,
\end{align}
with $\mathscr{M}(1)>0$ by $(M_1)$.
Fix $u\in C_c^\infty(\mathbb{R}^N,\mathbb{C})$, with $[u]_{s,A}=1$.
By~\eqref{eq3.8} and~\eqref{eqbu3.8} as $t\rightarrow\infty$
\begin{align*}
\mathcal{I}(tu)
&=\frac{1}{2}\mathscr{M}([tu]_{s,A}^2)+\frac{1}{2}\|tu\|_{L^2,V}^2-
\int_{\mathbb{R}^N}F(x,t|u|)dx\\
&\leq \frac{\mathscr{M}(1)}{2}t^{2\theta}[u]_{s,A}^{2\theta}
+\frac{1}{2}\|tu\|_{L^2,V}^2
-c_F t^\mu \|u\|_{L^\mu(\mathbb{R}^N)}^\mu
+\frac{M_1}{V_0} t^2\|u\|_{L^2,V}^2\\
&\leq \frac{\mathscr{M}(1)}{2}t^{2\theta}-c_FC_\mu^\mu t^\mu\|u\|_{s,A}^\mu
+\left(\frac{M_1}{V_0} +\frac{1}{2}\right) t^2\|u\|_{L^2,V}^2\\
&\leq\frac{\mathscr{M}(1)}{2}t^{2\theta}-c_FC_\mu^\mu t^\mu[u]_{s,A}^\mu
+\left(\frac{M_1}{V_0} +\frac{1}{2}\right) t^2\|u\|_{L^2,V}^2\\
&=\frac{\mathscr{M}(1)}{2}t^{2\theta}-c_FC_\mu^\mu t^\mu
+\left(\frac{M_1}{V_0} +\frac{1}{2}\right) t^2\|u\|_{L^2,V}^2\to-\infty,
\end{align*}
since $2\le2\theta<\mu$.
The assertion follows at once, taking $e=T_0u$, with $T_0>0$ large
enough.
\end{proof}

\section{Proof of Theorems~\ref{th1} and~\ref{th2}}\label{sec5}

The following standard Mountain Pass Theorem will
be used to get our main result.

\begin{theorem}\label{th4.1}
Let $J$ be a functional on a real Banach space $E$ and
of class $C^1(E,\mathbb{R})$. Let us assume that there exists $\alpha$, $\rho>0$ such that\\
$(i)$ $J(u)\geq\alpha$ for all $u\in E$ with $\|u\|=\rho$,\\
$(ii)$ $J(0)=0$ and $J(e)<\alpha$ for some $e\in E$ with $\|e\|>\rho$.\\
Let us define $\Gamma=\{\gamma\in C([0,1];E):\gamma(0)=0,\gamma(1)=e\}$, and
\begin{align*}
c=\inf_{\gamma\in \Gamma}\max_{t\in[0,1]}J(\gamma(t)).
\end{align*}
Then there exists a sequence $(u_n)_n$ in $E$ such
that $J(u_n)\rightarrow c$ and $J^\prime(u_n)\rightarrow 0$
in $E^\prime$, the dual space of $E$, as $n\to\infty$.
\end{theorem}

\subsection{Proof of Theorem \ref{th1}}
Taking into account  Lemmas \ref{lemma3.2} and \ref{lemma3.3},  by Theorem \ref{th4.1} there exists a sequence
$(u_n)_{n}\subset \Ma(\mathbb{R}^N,\mathbb{C})$ such that
$\mathcal{I}(u_n)\rightarrow c>0$ and $\mathcal{I}^\prime(u_n)\rightarrow 0$ as $n\to\infty$.
Then, in view of Lemma \ref{lemma3.1}, there exists a nontrivial
critical point $u_0\in \Ma(\mathbb{R}^N,\mathbb{C})$ of $\mathcal{I}$ with $\mathcal{I}(u_0)=c>0=\mathcal{I}(0)$.

Set $\mathscr{N}=\{u\in
\Ma(\mathbb{R}^N,\mathbb{C})\setminus\{0\}:\mathcal{I}^\prime(u)=0\}$.
Then $u_0\in \mathscr{N}\neq\emptyset$.
Next we show that $\mathcal{I}$ is coercive and bounded from below on $\mathscr{N}.$ Indeed, by $\mathcal{I}^\prime(u)=0$ and $(f_2)$, we get
\begin{align}\label{eq3.10}
\int_{\mathbb{R}^N}F(x,|u|)dx\leq\frac{1}{\mu}\int_{\mathbb{R}^N}f(x,|u|)|u|^2dx
=\frac{1}{\mu}\big(M([u]_{s,A}^2)[u]_{s,A}^2+\|u\|_{L^2,V}^2\big).
\end{align}
By using  \eqref{eq3.10}, $(M_2)$ and the fact that $2\le2\theta<\mu$,
for all $u\in\mathscr{N}$, we have
\begin{align*}
\mathcal{I}(u)&
\geq \frac{1}{2}\mathscr{M}(\|u\|_{s,A}^2)+\frac{1}{2}\|u\|_{L^2,V}^2
-\frac{1}{\mu}(M([u]_{s,A}^2)[u]_{s,A}^2+\|u\|_{L^2,V}^2)\\
&=\left(\frac{1}{2\theta}
-\frac{1}{\mu}\right)M([u]_{s,A}^2)[u]_{s,A}^2+\left(\frac{1}{2}
-\frac{1}{\mu}\right)\|u\|_{L^2,V}^2\ge0.
\end{align*}
Hence, by $(M_1)$ and $(M_3)$ for $u\in\mathscr{N}$
\begin{align}\label{eq3.11}
\mathcal{I}(u)&\ge\left(\frac{1}{2\theta}-\frac{1}{\mu}\right)
\cdot\left(\|u\|_{L^2,V}^2
+\begin{cases}\kappa[u]_{s,A}^2,&\mbox{if }[u]_{s,A}\ge1\\
m_0[u]_{s,A}^{2\theta},&\mbox{if }[u]_{s,A}\le1\end{cases}\right),
\end{align}
where $\kappa=\kappa(1)>0$ by $(M_1)$.
Hence in all cases, for all $u\in\mathscr{N}$
$$\mathcal{I}(u)\ge\min\{\kappa,m_0\}
\left(\frac{1}{2\theta}-\frac{1}{\mu}\right)\|u\|_{s,A}^2-1,$$
by the elementary inequality $t^\theta\ge t-1$ for all
$t\in\mathbb R^+_0$.
In particular, $\mathcal{I}$ is coercive and bounded from
below on $\mathscr{N}$.

Define $c_{\rm min}=\inf\{\mathcal{I}(u):u\in\mathscr{N}\}$.
Clearly, $0\le c_{\rm min}\leq \mathcal{I}(u_0)=c$.
Let $(u_n)_{n}$ be a minimizing
for $c_{\rm min},$ namely $\mathcal{I}(u_n)\rightarrow c_{\rm min}$ and $\langle \mathcal{I}^\prime(u_n),u_n\rangle=0$. Then,
since $\mathscr{N}$ is a complete metric space, by Ekeland's variational principle we can find a new minimizing sequence,
still denoted by $(u_n)_{n}$, which is a $(PS)$ sequence for
$\mathcal{I}$ at the level $ c_{\rm min}$.
Moreover, Lemma~\ref{lemma3.1} implies that $(u_n)_{n}$ has a
convergence subsequence, which we
still denote by $(u_n)_{n}$, such that
$u_n\rightarrow u$ in $\Ma(\mathbb{R}^N,\mathbb{C})$. Thus
$c_{\rm min}=\mathcal{I}(u)$ and
$\langle\mathcal{I}'(u),u\rangle=0$.

We claim that $c_{\rm min}>0$.
Otherwise, there is $(u_n)_{n}\subset \Ma(\mathbb{R}^N,\mathbb{C})\setminus\{0\}$ with $\mathcal{I}^\prime(u_n)=0$
and $\mathcal{I}(u_n)\rightarrow 0$. This via \eqref{eq3.11} implies that $\|u_n\|_{s,A}\rightarrow0$.
On the other hand, by \eqref{eqbu}, we have for any $\varepsilon\in (0,V_0)$
 \begin{equation*}
M([u_n]_{s,A}^2)[u_n]_{s,A}^2+\|u_n\|_{L^2,V}^2 =\int_{\mathbb{R}^N}f(x,|u_n|)|u_n|^2dx
 \leq \frac{\varepsilon}{V_0}\|u_n\|^2_{L^2,V}+C_\varepsilon C_p^p\|u_n\|_{s,A}^p.
 \end{equation*}
Thus, $M([u_n]_{s,A}^2)[u_n]_{s,A}^2+\Big(1-\varepsilon/V_0\Big)
\|u_n\|_{L^2,V}^2
\leq C_\varepsilon C_p^p\|u_n\|_{s,A}^p$.
Now take $N_1$ so large that $\|u_n\|_{s,A}\le1$ for all $n\ge N_1$.
Hence, $(M_3)$ implies that for all $n\ge N_1$
\begin{align*}
m_0[u_n]_{s,A}^{2\theta}+\big(1-\varepsilon/V_0\big)
\|u_n\|_{L^2,V}^{2\theta}
&\leq C_\varepsilon C_p^p\|u_n\|_{s,A}^p,
\end{align*}
that is
\begin{align*}
\min\big\{m_0,\big(1-\varepsilon/ V_0\big)\big\}\leq
 C_\varepsilon C_p^p\|u_n\|_{s,A}^{p-2\theta}.
\end{align*}
This is a contradiction since $2\theta<p$ and proves the claim.

Thus, $u$ is a nontrivial critical point of $\mathcal{I}$,
with $\mathcal{I}(u)=c_{\rm min}>0$.
Therefore, $u$ is a ground state solution of \eqref{eq1}.\qed

\subsection{Proof of Theorem \ref{th2}}
By $(f_5)$, $(V_1)$ and the H\"{o}lder inequality, for all
$u\in H_{s,A}^s(\mathbb{R}^N,\mathbb{C})$ we have
\begin{align*}
\mathcal{I}(u) &\geq \frac{1}{2}\mathscr{M}([u]_{s,A}^2)+\frac{1}{2}\|u\|_{L^2,V}^2-
\int_{\mathbb{R}^N}a(x)|u|^{q}dx\\
&\geq \frac{1}{2}\mathscr{M}([u]_{s,A}^2)
+\frac{1}{2}\|u\|_{L^2,V}^2-\|a\|_{L^{\frac{2}{2-q}}}
\|u\|_{L^2}^q\\
&\geq \frac{1}{2}\mathscr{M}([u]_{s,A}^2)+\frac{1}{4}\|u\|_{L^2,V}^2
+\frac{V_0}{4}\|u\|_{L^2}^2-\|a\|_{L^{\frac{2}{2-q}}}\|u\|_{L^2}^q\\
&\geq \frac{1}{2}\mathscr{M}([u]_{s,A}^2)+\frac{1}{4}\|u\|_{L^2,V}^2-C_0,\\
&C_0=\frac{\|a\|_{L^{\frac{2}{2-q}}}}{2q}(2q-1)
\left(\frac{2\|a\|_{L^{\frac{2}{2-q}}}}
{qV_0}\right)^{q/(2-q)}.
\end{align*}
As shown in \eqref{eq3.11}, this, $(M_1)$ and $(M_3)$
imply at once that for all $u\in H_{s,A}^s(\mathbb{R}^N,\mathbb{C})$
$$\mathcal{I}(u) \ge\frac{\min\{\kappa,m_0\}}4\,
\|u\|_{s,A}^2-1-C_0,$$
$\kappa=\kappa(1)$.
Hence $\mathcal{I}$ is coercive and
bounded below on $H_{s,A}^s(\mathbb{R}^N,\mathbb{C})$. Set
\begin{align*}
J(u)=\frac{1}{2}\mathscr{M}([u]_{s,A}^2)+\frac{1}{2}\|u\|_{L^2,V}^2,\quad
H(u)=\int_{\mathbb{R}^N}F(x,|u|)dx
\end{align*}
for  all $u\in H_{A,V}^s(\mathbb{R}^N)$.
Then $J$ is weakly lower semi--continuous in
$H_{A,V}^s(\mathbb{R}^N)$, since $\mathscr{M}$ is continuous and
monotone non--decreasing in $\mathbb R^+_0$. Moreover, by using a
similar discussion as \cite[Lemma~2.3]{PXZ2}, one can show that $H$ is weakly
continuous on $H_{A,V}^s(\mathbb{R}^N)$ under condition $(f_5)$.
Thus, $\mathcal{I}(u)=J(u)-H(u)$ is weakly lower semi--continuous
in $H_{A,V}^s(\mathbb{R}^N)$.
Then there exists $u_0\in H_{A,V}^s(\mathbb{R}^N)$ such that
$$\mathcal{I}(u_0)=\inf\{\mathcal{I}(u):u\in H_{A,V}^s(\mathbb{R}^N)\}.$$
Next we show $u_0\neq0$. Let $x_0\in \Omega$ and let $R>0$ such that $B_R(x_0)\subset \Omega$. Fix $\varphi\in C_0^\infty(B_{R}(x_0))$ with
$0\leq\varphi\leq 1$, $\| \varphi\|_{s,A}\leq C(R)$ and $\|\varphi\|_{L^q(B_R(x_0))}\neq 0$.
Then, by $(f_6)$ for all $t\in(0,\delta)$
\begin{align*}
\mathcal{I}(t\varphi)&\leq \frac{t^2}{2}\left(\sup_{0\leq \xi\leq (\delta C(R))^2}M(\xi)\right)[\varphi]_{s,A}^2+\frac{t^2}{2}\|\varphi\|_{L^2,V}^2
-t^{q}\int_{B_{R}(x_0)}a_0|\varphi|^qdx\\
&\leq \frac{t^2}{2}\left(\sup_{0\leq \xi\leq (\delta C(R))^2}M(\xi)+1\right)\|\varphi\|_{s,A}^2-
t^{q}a_0\|\varphi\|_{L^q(B_R(x_0))}.
\end{align*}
Since $1<q<2$, we get $\mathcal{I}(\bar t\varphi)<0$ by taking $\bar t>0$ small enough.
Hence $\mathcal{I}(u_0)\leq \mathcal{I}(\bar t\varphi)<0$,
and so $u_0$ is a nontrivial critical point.
In other words, $u_0$ is a nontrivial solution of \eqref{eq1}. \qed

\section{Proof of Theorem \ref{th3}}\label{sec6}

We first recall the following symmetric mountain pass theorem
in~\cite{RK}.

\begin{theorem}\label{th4}
Let $X$ be an infinite dimensional real Banach space.
Suppose that $J$ is in $C^1(X,\mathbb{R})$ and satisfies the following condition:\\
$(a)$ $J$ is even, bounded from below, $J(0)=0$ and $J$ satisfies
 the $(PS)$ condition;\\
$(b)$ For each $k\in\mathbb{N}$ there exists
$E_k\subset \Gamma_k$ such that
$\sup_{u\in E_k}J(u)<0$, where
\begin{align*}
\Gamma_k=\{E: E\mbox{ is closed symmetric subset
 of $X$  and } 0\notin E, \ \gamma(E)\geq k\}
\end{align*}
and $\gamma(E)$ is a genus of a closed symmetric set $E$.
Then $J$ admits a sequence of critical points $(u_k)_k$ such
that $J(u_k)\leq 0$, $u_k\neq0$ and $\|u_k\|\rightarrow0$
as $k\to\infty$.
\end{theorem}

Let $h\in C^1(\mathbb R^+_0,\mathbb{R})$ be a radial decreasing function such that $0\leq h(t)\leq 1$ for all
$t\in\mathbb R^+_0$, $h(t)=1$ for $0\leq t\leq 1$ and $h(t)=0$ for $t\geq 2$.
Let $\phi(u)=h(\|u\|^2_{s,A})$.
Following the idea of \cite{FV}, we consider the truncation functional
\begin{align*}
\mathcal{I}(u)=\frac{1}{2}\mathscr{M}([u]_{s,A}^2)
+\frac{1}{2}\|u\|_{L^2,V}^2-\phi(u)\int_{\mathbb{R}^N}F(x,|u|)dx.
\end{align*}
Clearly, $\mathcal{I}\in C^1(H_{A,V}^s(\mathbb{R}^N,\mathbb{C}),\mathbb{R})$ and
\begin{align*}
\langle \mathcal{I}^\prime(u),v\rangle&=M([u]_{s,A}^2)\langle L(u), v\rangle
+\Re\int_{\mathbb{R}^N}V(x)u\overline{v}dx\\
&\qquad-2\phi^\prime(u)\int_{\mathbb{R}^N}F(x,|u|)dx\cdot \langle L(u), v\rangle
-\phi(u)\Re\int_{\mathbb{R}^N}f(x,|u|)u\overline{v}dx
\end{align*}
for all $u,v\in H_{A,V}^s(\mathbb{R}^N,\mathbb{C})$. Here $L(u)$ is the linear functional
on $H_{A,V}^s(\mathbb{R}^N,\mathbb{C})$, introduced in~\eqref{Lu}.

\subsection{Proof of Theorem \ref{th3}}
For all $u\in H_{A,V}^s(\mathbb{R}^N,\mathbb{C})$, with $\|u\|_{s,A}\geq 2$, we get

$$\mathcal{I}(u)\geq \frac{1}{2}\mathscr{M}([u]_{s,A}^2)+\frac{1}{2}\|u\|_{L^2,V}^2\ge\frac{1}{2}\min\{\kappa,\,m_0\}\|u\|_{s,A}^2,$$
by  $(M_1)$ and $(M_3)$, where $\kappa=\kappa(1)$, as in the proof of
Theorem~\ref{th2}. Hence
$\mathcal{I}(u)\rightarrow \infty$ as
$\|u\|_{s,A}\rightarrow \infty$ and $\mathcal{I}$ is coercive and bounded from below on $H_{A,V}^s(\mathbb{R}^N,\mathbb{C})$.

Let $(u_n)_n$ be a $(PS)$ sequence, i.e. $\mathcal{I}(u_n)$ is bounded and $\mathcal{I}^\prime(u_n)\rightarrow 0$ as $n\rightarrow\infty$. Then the coercivity of $\mathcal{I}$
implies that $(u_n)_n$ is bounded in $H_{A,V}^s(\mathbb{R}^N,\mathbb{C})$. Without loss of generality, we assume that $u_n\rightharpoonup u$ in $H_{A,V}^s(\mathbb{R}^N,\mathbb{C})$ and
$u_n\rightarrow u$ a.e. in $\mathbb{R}^N$.
We now claim that
\begin{align}\label{bueq4}
\lim_{n\rightarrow\infty}
\int_{\mathbb{R}^N}(f(x,|u_n|)u_n-f(x,|u|)u)(\overline{u_n-u})dx=0.
\end{align}
Clearly, $|f(x,t)t|\leq C(|t|+|t|^{p-1})$ for  all
$x\in\mathbb{R}^N$ and $t\in\mathbb{R}^+_0$ by $(f_1)$.
Using the H\"{o}lder inequality, we obtain
\begin{align}\label{eqbu3.3}
&\int_{\mathbb{R}^N}|(f(x,|u_n|)u_n-f(x,|u|)u)(\overline{u_n-u})|dx\nonumber\\
&\leq \int_{\mathbb{R}^N}
C[|u_n|+|u|+|u_n|^{p-1}+|u|^{p-1}]|u_n-u|dx\\
&\leq C(\|u_n\|_{L^2}+\|u\|_{L^2})\|u_n-u\|_{L^2}
+C(\|u_n\|_{L^p}^{p-1}+\|u\|_{L^p(\mathbb{R}^N)}^{p-1})\|u_n-u\|_{L^p(\mathbb{R}^N)}
\nonumber\\
&\leq C(\|u_n-u\|_{L^2}+\|u_n-u\|_{L^p}).\nonumber
\end{align}
Lemma \ref{lemma2.3} guarantees that $|u_n|\rightarrow |u|$ in $L^p(\mathbb{R}^N)$  and $|u_n|\rightarrow |u|$ in $L^2(\mathbb{R}^N)$.
Hence,  $u_n\to u$ in $L^p(\mathbb{R}^N,\mathbb{C})$ and in $L^2(\mathbb{R}^N,\mathbb{C})$ by the Brezis--Lieb lemma.
Inserting these facts in \eqref{eqbu3.3}, we get the desired claim \eqref{bueq4}.

Now, $\langle \mathcal{I}^\prime(u_n)-\mathcal{I}^\prime(u),u_n-u\rangle\rightarrow0$,
since $\mathcal{I}^\prime(u_n)\rightarrow0$ and $u_n\rightharpoonup u$ in $H_{A,V}^s(\mathbb{R}^N,\mathbb{C})$.
By \eqref{bueq4}, we have as $n\to\infty$
\begin{align*}
o(1)&=\langle \mathcal{I}^\prime(u_n)-\mathcal{I}^\prime(u),u_n-u\rangle=M([u_n]_{s,A}^2)\langle L(u_n), u_n-u\rangle-M([u]_{s,A}^2)\langle L(u), u_n-u\rangle\\
&+\Re\int_{\mathbb{R}^N}V(x)(u_n-u)\overline{(u_n-u)}dx
-2\phi^\prime(u_n)\int_{\mathbb{R}^N}F(x,|u_n|)dx\cdot \langle L(u_n), u_n-u\rangle\\
&-2\phi^\prime(u)\int_{\mathbb{R}^N}F(x,|u|)dx\cdot \langle L(u), u_n-u\rangle
-\phi(u_n)\Re\int_{\mathbb{R}^N}f(x,|u_n|)u_n(\overline{u_n-u})dx\\
&-\phi(u)\Re\int_{\mathbb{R}^N}f(x,|u|)u(\overline{u_n-u})dx.
\end{align*}
From  $(f_7)$  and the facts that $u_n\rightharpoonup u$
in $H_{A,V}^s(\mathbb{R}^N,\mathbb{C})$ and $\phi^\prime\leq 0$ it
follows that
\begin{align}\label{PPP}
0\le M([u_n]_{s,A}^2)\langle L(u_n)-L(u), u_n-u\rangle
+\Re\int_{\mathbb{R}^N}V(x)(u_n-u)\overline{(u_n-u)}dx\leq o(1).
\end{align}
We divide the proof into two parts.\smallskip

\noindent {\em Case $\inf_{n\in\mathbb N}[u_n]_{s,A}=d>0$.}
By $(M_1)$, there exists $\kappa=\kappa(d)>0$ with $M(t)\geq \kappa>0$ for all $t\geq d$.
This, together with \eqref{PPP}, implies that
\begin{align*}
\lim_{n\rightarrow\infty}\bigg[\iint_{\mathbb{R}^{2N}}\frac{|u_n(x)-u(x)-e^{{\rm i}(x-y)\cdot A(\frac{x+y}{2})}(u_n(y)-u(y))|^2}{|x-y|^{N+2s}}dxdy
+\int_{\mathbb{R}^N}V(x)|u_n-u|^2dx\bigg]=0.
\end{align*}
Hence $u_n\rightarrow u$ in $H_{A,V}^s(\mathbb{R}^N,\mathbb{C})$.
\smallskip

\noindent{\em Case $\inf_{n\in\mathbb N}[u_n]_{s,A}=0$}.
If $0$ is an isolated point for $([u_n]_{s,A})_n$, then there is a
subsequence $([u_{n_k}]_{s,A})_k$ such that
$\inf_{k\in\mathbb N}[u_{n_k}]_{s,A}=d>0$ and one can proceed as before.

If, instead, $0$ is an accumulation point for $([u_n]_{s,A})_n$, there is a subsequence,
still labeled as $(u_n)_n$, such that
$[u_n]_{s,A}\to 0$ and $u_n\to 0$ in $L^{2^*_s}(\mathbb{R}^N)$ as $n\to\infty$ and
again \eqref{PPP} implies at once that
$u_n\rightarrow 0$ in $H_{A,V}^s(\mathbb{R}^N,\mathbb{C})$, since $\langle L(u_n)-L(u), u_n-u\rangle\to0$ and $M([u_n]_{s,A}^2)\to M(0)\ge0$.

In conclusion, $\mathcal{I}$ satisfies the $(PS)$ condition in $H_{A,V}^s(\mathbb{R}^N,\mathbb{C})$.
For each $k\in\mathbb N$, we take $k$ disjoint open sets $K_i$
such that $\bigcup_{i=1}^k K_i\subset\Omega$. For each
$i=1,\cdots,k$ let $u_i\in (H_{A,V}^s(\mathbb{R}^N,\mathbb{C})
\bigcap C_0^\infty(K_i,\mathbb{C}))\backslash\{0\}$, with $\|u_i\|_{s,A}=1$,
and $W_k={\rm span}\{u_1,u_2,\cdots,u_k\}.$
Therefore, for any $u\in W_k$, with
$\|u\|_{s,A}=\rho\leq 1$ small enough, we obtain  by $(f_7)$,
being $q\in(1,2)$,
\begin{align*}
\mathcal{I}(u)&\leq \frac{1}{2}\left(\max_{0\leq t\leq 1}M(t)\right)[u]_{s,A}^2+\frac{1}{2}\|u\|_{L^2,V}^2
-\int_{\Omega}a_1|u|^qdx\\
&\leq \frac{1}{2}
\left(1+\max_{0\leq t\leq 1}M(t)\right)\|u\|_{s,A}^2
-C_k^qa_1\|u\|_{s,A}^q\\
&=\frac{1}{2}\left(1+\max_{0\leq t\leq 1}M(t)\right)\rho^2-
C_k^q a_1\rho^q<0,
\end{align*}
where $C_k>0$ is a constant such that $\|u\|_{L^q(\mathbb{R}^N,\mathbb{C})}
\le C_k\|u\|_{s,A}$ for all $u\in W_k$, since all norms on~$W_k$ are equivalent.
Therefore, we deduce
\begin{align*}
\{u\in W_k:\|u\|_{s,A}=\rho\}\subset \{u\in W_k:\mathcal{I}(u)<0\}.
\end{align*}
Obviously, $\gamma(\{u\in W_k:\|u\|_{s,A}=\rho\})=k$, see~\cite{chang}.
Hence by the monotonicity of the
genus $\gamma$, cf. \cite{krasnoselskii}, we obtain
\begin{align*}
\gamma(u\in W_k:\mathcal{I}(u)<0)\geq k.
\end{align*}
Choosing $E_k=\{u\in W_k:\mathcal{I}(u)<0\}$, we have $E_k\subset \Gamma_k$ and
$\sup_{u\in\Gamma_k}\mathcal{I}(u)<0$. Thus, all the assumptions of Theorem~\ref{th4} are satisfied,
Hence, there exists a sequence $(u_k)_k$ such that
$$\mathcal{I}(u_k)\leq 0,\quad \mathcal{I}^\prime(u_k)=0,\quad\mbox{and}\quad
\|u_k\|_{s,A}\rightarrow 0\mbox{ as }k\rightarrow\infty.$$
Therefore, we can take $k$ so large  that
$\|u_k\|_{s,A}\leq 1$, and so these infinitely many functions $u_k$
are solutions of~\eqref{eq1}.\qed

\bigskip
\bigskip

\medskip

\begin{thebibliography}{99}
\bibitem{r8}
D. Applebaum,
L\'{e}vy processes -- From probability to finance quantum groups,
{\em Notices Amer. Math. Soc.} {\bf 51} (2004), 1336--1347.

\bibitem{GAAS}
G.\ Arioli, A.\ Szulkin,
A semilinear Schr\"{o}dinger equation in the presence of a magnetic field,
{\em Arch. Rational Mech. Anal.} {\bf 170} (2003), 277--295.

\bibitem{AFP} G. Autuori, A. Fiscella, P. Pucci, Stationary Kirchhoff problems
involving a fractional
operator and a critical nonlinearity, {\em Nonlinear Anal.} {\bf125} (2015), 699--714.

\bibitem{BW}
T.\ Bartsch, Z.-Q.\ Wang, Existence and multiplicity results for some superlinear elliptic problems on $\mathbb{R}^{N}$.
{\em Comm. Partial Differential Equations} {\bf 20} (1995), 1725--1741


\bibitem{r3}
L. Caffarelli,
Nonlocal equations, drifts and games,
{\em Nonlinear Partial Differential Equations}, Abel Symposia {\bf 7} (2012), 37--52.


\bibitem{r7}
L. Caffarelli, L. Silvestre,
An extension problem related to the fractional Laplacian,
{\em Comm. Partial Differential Equations} {\bf 32} (2007), 1245--1260.

\bibitem{capu}
M. Caponi, P. Pucci,
Existence theorems for entire solutions of stationary Kirchhoff fractional $p$--Laplacian equations,
{\em Ann. Mat. Pura Appl.} DOI 10.1007/s10231-016-0555-x.

\bibitem{chang}
K.C. Chang, Critical point theory and applications, Shanghai Scientific and Technology Press, Shanghai, 1986.

\bibitem{SCSS}
S. Cingolani, S. Secchi,
Semiclassical limit for nonlinear Schrodinger equations with
electromagnetic fileds, {\em J. Math. Anal. Appl.} {\bf 275} (2002), 108--130.

\bibitem{colpuc}
F. Colasuonno, P. Pucci, Multiplicity of solutions
for $p(x)$--polyharmonic Kirchhoff equations,
{\em Nonlinear Anal.} 74 (2011), 5962--5974.

\bibitem{dAS}
 P. D'Ancona, S. Spagnolo,
Global solvability for the degenerate Kirchhoff equation with real
analytic data, {\em Invent. Math.} {\bf 108} (1992), 247--262.

\bibitem{PDMS}
P. d'Avenia, M. Squassina,
Ground states for fractional magnetic operators, preprint, (2016)
\href{http://arxiv.org/abs/1601.04230}{arXiv link}

\bibitem{JDJV}
J. Di Cosmo, J. Van Schaftingen,
Semiclassical stationary states for nonlinear
Schr\"{o}diner equations under a strong extenal magnetic field,
{\em J. Differential Equations} {\bf 259} (2015), 596--627.

 \bibitem{r28}
 E. Di Nezza, G. Palatucci, E. Valdinoci,
 Hitchhiker guide to the fractional Sobolev spaces, {\em Bull. Sci. Math.} {\bf 136} (2012), 521--573.

 \bibitem{DMV}
 S. Dipierro, M. Medina, E. Valdinoci,
Fractional elliptic problems
with critical growth in the whole of $\mathbb R^n$, arXiv:1506.01748, pp. 122.

 \bibitem{r15}
 S. Dipierro, G. Palatucci, E. Valdinoci,
 Existence and symmetry results for a  Schr\"{o}dinger type problem involving the fractional Laplacian,
 {\em Le Matematiche} {\bf 68} (2013), 201--216.


\bibitem{r14}
P. Felmer, A. Quaas, J. Tan,
Positive solutions of the nonlinear Schr\"{o}dinger equation with the fractional Laplacian,
{\em Proc. Roy. Soc. Edinburgh Sect. A} {\bf 142} (2012), 1237--1262.

\bibitem{FV}
A. Fiscella, E. Valdinoci, A critical Kirchhoff type
problem involving a nonlocal operator, {\em Nonlinear Anal.} {\bf 94}
(2014), 156--170.

\bibitem{ichi1}
T.\ Ichinose, Essential selfadjointness of the Weyl quantized relativistic Hamiltonian, {\em Ann. Inst. H. Poincar\'e Phys. Th\'eor.} {\bf 51} (1989), 265--297.

\bibitem{ichi2}
T.\ Ichinose, Magnetic relativistic Schr\"odinger operators and imaginary-time path integrals, Mathematical physics, spectral theory and stochastic analysis, 247--297, Oper. Theory Adv. Appl. {\bf 232}, Birkh\"auser/Springer Basel AG, Basel, 2013.

\bibitem{ichi3}
T.\ Ichinose, H.\ Tamura, Imaginary-time path integral for a relativistic spinless particle in an electromagnetic field, {\em Comm. Math. Phys.} {\bf 105} (1986), 239--257.

\bibitem{purice}
V.\ Iftimie, M.\ M\u{a}ntoiu, R.\ Purice, Magnetic pseudodifferential operators, {Publ. Res. Inst. Math. Sci.} {\bf 43} (2007), 585--623.

\bibitem{RK}
R. Kajikiya,
A critical point theorem related to the symmetric mountain pass lemma
and its applications to elliptic equations,
{\em J. Funct. Anal.} {\bf 225} (2005), 352--370.

\bibitem{krasnoselskii}
M.A. Krasnoselskii,
Topological methods in the theory of nonlinear integral equations,
MacMillan, New York, 1964.

\bibitem{KK}
K. Kurata,
Existence and semi-classical limit of the least energy solution to a nonlinear
Schr\"{o}dinger equation  with electromagnetic fileds, {\em Nonlinear Anal.} {\bf 41} (2000), 763--778.


\bibitem{laskin1}
N. Laskin,
Fractional quantum mechanics and L\'{e}vy path integrals,
{\em Phys. Lett. A} {\bf 268} (2000), 298--305.

\bibitem{laskin2}
N. Laskin,
Fractional Schr\"{o}dinger equation, {\em  Phys. Rev. E}
{\bf 66} (2002), 056108.

\bibitem{Lions}
P.L. Lions,
Sym\'{e}trie et compacit\'{e} dans les espaces de Sobolev, {\em
J. Funct. Anal.} {\bf 49} (1982), 315--334.

\bibitem{XMTZ}
X. Mingqi, G. Molica Bisci, G. Tian, B. Zhang,
Infinitely many solutions for the stationary Kirchhoff problems
involving the fractional $p$--Laplacian, {\em Nonlinearity} {\bf 29} (2016), 357--374.

\bibitem{MBRS}
G. Molica Bisci, V.D. Radulescu, R. Servadei, {\em Variational methods for nonlocal
fractional equations},
Encyclopedia of Mathematics and its Applications, {\bf 162},
Cambridge University Press, Cambridge, 2016.

\bibitem{PS}
P. Pucci, S. Saldi, Critical stationary Kirchhoff equations in $\mathbb{R}^N$ involving nonlocal operators, {\em Rev. Mat. Iberoam.}
{\bf 32} (2016), 1--22.

\bibitem{PXZ}
P. Pucci, M. Xiang, B. Zhang,
Multiple solutions for nonhomogeneous Schr\"{o}dinger-Kirchhoff type equations
involving the fractional $p$-Laplacian in ${\mathbb {R}}^N$,
{\em Calc. Var. Partial Differential Equations} {\bf 54} (2015), 2785--2806.

 \bibitem{PXZ2}
 P. Pucci, M. Xiang, B. Zhang,
 Existence and multiplicity of entire solutions for fractional $p$-Kirchhoff equations,
 {\em Adv. Nonlinear Anal.} {\bf 5} (2016), 27--55.

 \bibitem{r17}
S. Secchi,
Ground state solutions for the fractional Schr\"{o}dinger in
$\mathbb{R}^N$, {\em J. Math. Phys.} {\bf 54}, 031501, 2013

\bibitem{MS}
M. Squassina,
Soliton dynamics for the nonlinear Schr\"{o}dinger equation
with magnetic field,
{\em Manuscripta Math.} {\bf 130} (2009), 461--494.

\bibitem{MSBV}
M. Squassina, B. Volzone,
Bourgain-Brezis-Mironescu formula for magnetic operators,
{\em Comptes Rendus Mathematique}, DOI	10.1016/j.crma.2016.04.013.

\bibitem{XZF1}
M. Xiang, B. Zhang, M. Ferrara,
Existence of solutions for Kirchhoff type
problem involving the non-local fractional $p-$Laplacian, {\em J. Math. Anal. Appl.} {\bf 424} (2015), 1021--1041.

\bibitem{XZF2}
M. Xiang, B. Zhang, M. Ferrara,
Multiplicity results for the nonhomogeneous fractional $p$--Kirchhoff equations with concave--convex nonlinearities,
{\em Proc. A} {\bf 471} (2015), 2177, 14 pp.

\bibitem{XZG}
M. Xiang, B. Zhang, X. Guo, Infinitely many solutions for a fractional Kirchhoff type problem via Fountain Theorem,
{\em Nonlinear Anal.} {\bf 120} (2015), 299--313.
\end{thebibliography}
\end{document}